\documentclass{article}
\usepackage[utf8]{inputenc}
\usepackage[tbtags]{amsmath}
\usepackage{amsfonts,amssymb,mathrsfs,amscd}

\def \Z {{\mathbf {Z}}}

\def \B {{\cal B}}
\def \T {{\mathbf  T}}

\def \R {{\mathbf  R}}

\def\eps{\varepsilon}

\begin{document}

\title{\LARGE    \bf Multiple mixing, 75 years of  Rokhlin's problem}
\author{\LARGE   Valery \,V.~Ryzhikov}

\date{12.11.2024}

\maketitle


\Large

\begin{abstract}{In 1949, V.A. Rokhlin introduced  new  invariants  of measure-preserving transformations, called k-fold  mixing.  Does mixing imply k-fold mixing? -- this problem remains open. We recall shortly  some results and discuss related problems.}
\end{abstract}



\section{Introduction }
In  \cite{Ro} (1949), V.A. Rokhlin wrote:

\vspace{2mm} {\large
\it The proposed work arose from the author's attempts to solve a well-known
spectral problem in the theory of dynamical systems: are there
metrically different dynamical systems with the same continuous
(in particular, Lebesgue) spectrum? Using new metric invariants introduced
for this purpose, the author tried to find
among the ergodic automorphisms of compact commutative groups
automorphisms of different metric types. It turned out, however, that for all the indicated automorphisms the new invariants are absolutely identical.}

\vspace{2mm}
Rokhlin's  new invariants  are called mixing of multiplicity $k> 1$. An automorphism $T$ of a probability space $(X,\B,\mu)$ has $k$-fold mixing if for any measurable sets $A_0,A_1,\dots, A_k$ with $m_1,\dots m_k\, \to\,+\infty$ we have
$$\mu(A_0\cap T^{m_1}A_1\cap T^{m_1+m_2}A_2\cap \dots T^{m_1+\dots+m_k}A_k)\ \to \ \mu(A_0)\mu(A_1)\dots \mu(A_k).$$
\\ \\ 
{\ \ \ \ \  \LARGE  \bf  Rokhlin's problem: \it  does  1-mixing imply k-fold mixing?}
\newpage
In \cite{Ro} multiple mixing was proved for ergodic automorphisms of compact commutative groups  (they preserve the Haar probability measure). 
Rokhlin's problem can be naturally generalized to group actions. An elegant counterexample  for $\Z^2$-actions generated by automorphisms of compact commutative groups was given by F. Ledrappier\cite{L}. 

\vspace{3mm}\it
Rokhlin's idea to consider automorphisms of compact commutative groups 

for counterexamples is successful for $\Z^n$-actions, excluding $n=1$. 

\rm
\vspace{3mm}
The  first example of a pair of non-isomorphic automorphisms with the same spectrum was proposed by Anzai \cite{An}. The desired pair is the appropriate skew products $T$ and $T^{-1}$. As for Rokhlin's idea to use multiple asymptotic  properties to distinguish automorphisms with the same spectrum was right as well. Now we know that  for the  generic automorphism $T$ there exists a sequence $m_j\to\infty$ such that
$$4\mu(A\cap T^{m_j}A\cap T^{3m_j}A)\to \ \mu(A)+ \mu(A)^2 + 2\mu(A)^3,$$
and 
$$\mu(A\cap T^{-m_j}A\cap T^{-3m_j}A)\to \ \mu(A)^2$$
for any measurable set $A$ (see \cite{R23}).
\rm Automorphisms $T$ and $T^{-1}$ have  different 2-fold  invariants and the same spectrum. 

\vspace{3mm} 
The problem on multiple mixing is relevant even  in the class of weakly mixing automorphisms not necessary 1-mixing. \it

\vspace{4mm} 
Let for   all $A,B,C\in \B$ 
$$\mu(A\cap T^{m_i}B),\ \mu(A\cap T^{n_i}B),\ \mu(A\cap T^{n_i+m_i}B)\ \to\ \mu(A)\mu(B).$$
Is it true that 
$$\mu(A\cap T^{m_i} B\cap T^{m_i+n_i}C)\ \to \ \mu(A)\mu(B)\mu(C)?$$
\rm

\newpage
\bf A bit of philosophy. \rm The mixing property can be interpreted as the asymptotic independence of an event now and an event in the distant past.  If our random process is stationary, will the event now,  be asymptotically independent of a pair of events: a distant one and a very distant one? The question is open, and only under some  specific conditions we have  some   positive results. 
It is a bit surprising  that the slower the mixing, the more reliably we observe the multiple mixing.  The absence of multiple mixing for mixing commutative actions in fact requires  1-mixing to be  very good.

\vspace{4mm}
\bf Brief summary. \rm Mixing imply multiple  mixing for

-- automorphisms with singular spectrum \cite{Ho};

-- unipotent flows \cite{Ma}--\cite{St};

-- flows with the Ratner property \cite{Ra} or its analogues \cite{FK};

-- flows of positive local rank \cite{00}, quasi-simple flows \cite{RT};

-- automorphisms of finite rank \cite{Ka}, \cite{93}. \rm

\vspace{2mm}
The multiple mixing property is usually established by finding some algebraic, spectral, approximation properties  that force it. Our goal is  to discuss  such relationships.

\section{ Homoclinic group, commutation relations}
The Gaussian  actions are the  result of the embedding of the orthogonal group into the automorphism group of a space with Gaussian measure. The  Poisson suspensions are the  image of the embedding of infinite automorphism group into the group of    Poisson measure space  automorphism.
For Gaussian mixing automorphisms, multiple mixing was proved by V.P. Leonov \cite{Leo}. For Poisson mixing systems, it follows in fact  from the definition. 
We will not give definitions of Gaussian and Poisson suspensions, but we  say that they have the Gordin  ergodic group \cite{19}, which ensures the multiple mixing property. 

Let $T$ be an automorphism of a probability space, its homoclinic group $H(T)$ is 
$$H(T)=\{S: \ T^{-n}ST^n\to I,\ n\to\infty\}.$$
M.I. Gordin showed that the ergodicity of $H(T)$ implies $T$ be mixing.   In fact   mixing can be replaced by multiple mixing.

\vspace{2mm}
\bf Theorem 2.1. \it An automorphism with  ergodic homoclinic group is mixing of all orders.\rm

Proof. The group $H(T)$ is full (see \cite{19}).  If it is ergodic, then there is an ergodic automorphism  $S\in H(T)$. Fix an integer $i$. For any measurable sets $A,B,C$ 
$$
\mu(S^{-i}A\cap T^{m} B\cap T^{m+n}C)-
\mu(A\cap T^{m} B\cap T^{m+n}C)\to 0,\ m, n\to +\infty.$$
Indeed, from
$$\mu (T^{-m} S^iT^{m}B\Delta B)\to 0, \ \ \ T^{-m-n}S^iT^{m+n}C \Delta C)\to 0$$
we get 
$$\mu(A\cap T^{m}T^{-m} S^iT^{m} B\cap T^{m+n}T^{-m-n}S^iT^{m+n}C) -
\mu(A\cap T^{m} B\cap T^{m+n}C)\to 0.$$
  Then $$ \mu(A\cap T^{m} B\cap T^{m+n}C) -\frac 1 N \sum_{i=1}^N \mu(S^{-i}A\cap T^{m} B\cap T^{m+n}C)\to 0,\ m, n\to +\infty.$$

From the ergodicity of the automorphism $S$ we have
$$ \left\|\frac 1 N \sum_{i=1}^N S^{-i}\chi_A - \mu(A)\chi_X\right\| \to \, 0, N\to\infty,$$
and from  1-mixing  of $T$ we have 
$$ \mu(T^{m} B\cap T^{m+n}C)\to \mu(B)\mu(C).$$
 So for every $\eps>0$ there is $N$ such that for all sufficiently large $ m,n$ $$ \left|\mu(A)\mu(B)\mu(C) - \frac 1 N \sum_{i=1}^N \mu(S^{-i}A\cap T^{m} B\cap T^{m+n}C)\right| < \eps.$$
Thus, 
$$\mu(A\cap T^m B\cap T^{m+n}C)\to \mu(A)\mu(B)\mu(C), \ m,n\to +\infty.$$
Similarly, by induction, we get $T$ is $k$-fold mixing for all $k$.

\vspace{3mm}
\bf Problem. \rm  It would be interesting to find  non-isomorphic to Gaussian and Poisson suspensions automorphisms $T$ with    ergodic group $H(T)$. Probably,  King's  transformation from  \cite{KG} is so.

\vspace{3mm} 
{\bf Commutation relations.}
The following method works for some  flows with   Lebesgue spectrum 
(for example, for  unipotent flows). Multiple mixing of such flows is a consequence of nontrivial commutation relations, see  \cite{Ma}, \cite{91}, \cite{Mo}.  Let us consider   Heisenberg group actions.
\medskip

\bf Theorem 2.2. \it Let measure-preserving ergodic flows
$\Psi_r ,\Phi_s, T_t $ satisfy the relation
$\Psi_r\Phi_s= T_{rs}\Phi_s\Psi_r,$
where $T_c$ is a weakly mixing flow commuting with
$\Psi_r$ and $\Phi_s$. Then the flows $\Psi_r$ and $\Phi_s$ have mixing of all orders. \rm
\medskip

Proof. Define an operator $P$:
$$
\left< A,P(B\otimes C)\right> =
\lim_{i'\to\infty} \left< A,\,\Psi_{m(i')}B \Psi_{n(i')}C)\right>.
$$
Suppose that $ 0 < m(i') < n(i')$ and there exists 
$$\lim_{i} m(i)/n(i) = \alpha,$$
where  $i$ is a subsequence of $i'$.
Note that for $s\to +0$ we have
$$
\left< \Phi_sA,\,(\Psi_{m(i)}\Phi_sB)(\Psi_{n(i)}\Phi_sC)\right> \to
\left< A,\,(\Psi_{m(i)}B)(\Psi_{n(i)}C)\right>.
$$ Therefore, $$ \lim_{\eps\to 0}\lim_{i\to\infty} \frac{1}{\eps}\int_{0}^\varepsilon \left< \Phi_sA,\,(\Psi_{m(i)}\Phi_sB)(\Psi_{n(i)}\Phi_sC)\right> ds = $$ $$ \lim_{i\to\in fty} \left< A,\,(\Psi_{m(i)}B)(\Psi_{n(i)}C)\right>.
$$
Commutation relations give
$$
\left< \Phi_sA,\,(\Psi_{m(i)}\Phi_sB)(\Psi_{n(i)}\Phi_sC)\right> =
\left< A,\,(\Psi_{m(i)}T_{sm(i)}B)(\Psi_{n(i)}T_{sn(i)}C)\right>.
$$ Thus, $ \left< A,P(B\otimes C)\right> = $ $$ \lim_{\eps\to 0}\lim_{i\to\infty} \frac{1}{\eps}\int_{0}^\varepsilon \left< A,\Psi_{m(i)}T_{sm(i)}B \Psi_{n(i)}T_{sn(i)}C )\right>ds.
$$
The expression
$$ \frac{1}{\eps}\int_{0}^\varepsilon
\left< A,\Psi_{m(i)}T_{sm(i)}B
\Psi_{n(i)}T_{sn(i)}C)\right> ds
$$
for fixed $r$ and large $i$ differs little from
$$ \frac{1}{\eps}\int_{0}^\eps
\left< A,\Psi_{m(i)}T_{(s-\frac{r}{n(i)})m(i)}B
\Psi_{n(i)}T_{(s-\frac{r}{n(i)})n(i)}C)\right> ds.
$$
Therefore, for  $r\in\bf R$, we get
$$
\left< A,P(B\otimes C)\right> =
\lim_{\eps\to 0}\lim_{i\to\infty} \frac{1}{\eps}\int_{0}^\varepsilon
\left< A,\Psi_{m(i)}T_{sm(i)}T_{ar}B
\Psi_{n(i)}T_{sn(i)}T_{r}C)\right> ds.
$$
This leads to the equality
$
P(T_{ar}\otimes T_{r} )= P.
$
Since $T_{r}$ is weakly mixing, for $a\neq 0$ the product $T_{ar}\otimes T_{r}$ is ergodic, so  we get  $$ {\rm Im} \ P^{\ast}=\{Const\}.$$
If $a=0$, the last  is true as  well (exercise).  Thus,
$$\left< A,P(B\otimes C)\right> = \mu(A)\mu(B)\mu(C).
$$

\vspace{2mm}
\bf Problem. \rm Is it true that such mixing $T_{t}$ (the center of the Heisenberg group action)  is  mixing of all orders?

\section{Joinings and multiple mixing}
Some properties of actions formulated in terms of joinings (or Markov intertwining operators) imply  multiple mixing property. 

A \it joining \rm  of actions $\Psi_1,\dots,\Psi_n$ is a measure on $X^n =X_1\times \dots\times X_n$ ($X_i=X$) whose projections onto the edges of the cube $X^n$ are equal to $\mu$, and the measure is invariant under the diagonal action of the product $\Psi_1\times \dots\times \Psi_n$. If $\Psi_1,\dots,\Psi_n$ are copies of a single action, such a joining is called a self-joining.

We say that the action $\Psi$ belongs to the class $S(m, n)$, $n>m>1$ (or has the property $S(m,n)$), if every self-joining of order $n>2$ such that all projections onto the $m$-dimensional faces of the cube $X^n$ are equal to $\mu^m$, is trivial, i.e. coincides with the measure $\mu^n$,
the product of $n$ copies of the measure $\mu$. In what follows, we use the notation $S_n=S(n-1,n)$, $n>2$.

\bf Joinings and multiple mixing. \rm A useful observation, arrived at independently by several researchers, relates the properties of $S_n$ to multiple mixing:

\it if a mixing automorphism has mixing of multiplicity
$n-1$ and the property $S_{n+1}$, then it has multiple mixing
of order $n$. \rm

Let for any measurable $A,B,C$ we have
$$\mu(A\cap T^{m_i} B\cap T^{m_i+n_i}C)\to \nu(A\times B\times C), $$
Property $S_3$ of automorphism $T$ implies
$$\nu(A\times B\times C) =\mu(A)\mu(B)\mu(C).$$
Indeed, $\nu$ is a normalized measure on the semiring of cylinders
of the form $A\times B\times C$.
It is verified that
$$\nu(A\times B\times C) = \nu(TA\times TB\times TC),$$
$$\nu(A\times B\times X)=\mu(A)\mu(B), \ \nu(X\times B\times C) = \mu(B)\mu(C),$$
$$ \nu(A\times X\times C) = \mu(A)\mu(C).$$
Such measures are called self-joinings with pairwise independence.

If the mixing automorphism $T$ has the property $S_3$, then for any measurable sets $A,B,C\in \B$
$$\nu(A\times B\times C)=\mu(A)\mu(B)\mu(C),$$
hence, we obtain that $T$ has mixing of multiplicity 2. Similarly, we see that

\it  a mixing automorphism  with   the properties $S_3,\dots,S_{k+1}$ 
is $k$-fold mixing. \rm

\newpage
\bf del Junco-Rudolph problem \cite{JR}. \rm
Is there  a weakly mixing automorphism with zero entropy that does not have the property 
$PID=\bigcap_{p>2}\, S_p$? 

The following result, although proven in one line, was a pleasant surprise for experts.

\vspace{3mm}
\bf Theorem (J. King \cite{Ki}). \ \it $S_3$ and $S_4$ together  imply  $PID$.  \rm 

\vspace{3mm}
In fact the properties $S_{2m}$, $m>1$, for all group actions are equivalent to each other, are $PID$, and force all the properties $S_{2m-1}$, $m>1$ (see \cite{97}). 
 
 In \cite{RT}, \cite{93}, \cite{GHR},   \cite{97s}  for certain classes of actions it is shown that $S_3$ implies  $S_4$.
There are group actions that have all $S_{2m-1}$, but do not have  $PID$. 
Example: consider a compact commutative group $X=\Z_2^{\Z}$ with Haar measure $\mu$, on $(X,\mu)$ we define the action $\Psi$ generated by all shifts on the group $X$ and all automorphisms of the group $X$. It turns out that the action $\Psi$ and many of its finitely generated non-commutative subactions have all the properties of $S_{2m-1}$, but do not have the property $S_{4}$.

\vspace{2mm}
\bf Flows. \rm For   some important classes of flows the multiple mixing can be proved by use of Rantner property  \cite{Ra}, \cite{Th} or its analogues \cite{FK}. Here we are dealing with  quasi-simple flows,  possessing  $PID$, see \cite{RT}.  
The mixing flows of positive local rank are also quasi-simple, then  mixing of all orders
 \cite{00}.  Flows, unlike automorphisms, have the advantage of allowing small perturbations of joinings.  Therefore, some results for flows are better than for automorphisms.
One can safely conjecture (the author did this more than 30 years ago in one of his publications) that for flows, the Rokhlin problem has positive solution even in the class of $\R^n$-actions. This conjecture  is slightly confirmed by the work  \cite{Ti}, where S. Tikhonov shows that  mixing $\Z^n$-action with extreme deviation from 2-fold  mixing  has to have trivial  centralizer. So it cannot be any $\R^n$-flow. \rm

\newpage
\section{Joinings and spectrum}
\bf Theorem (B. Host \cite{Ho}).
\it Automorphisms with  singular spectra have  $PID$. \rm

\vspace{2mm}
\bf Corollary. \it Mixing automorphisms  with  singular spectra are mixing of all orders. \rm

\vspace{2mm}
In brief, the scheme of Host's reasoning is as follows (in the paper \cite{Ho} a more general situation of joining of three transformations $R,S,T$ is considered, but we will restrict ourselves to the special case when
$R=S=T$).
Let  $\tau$ be a measure on the torus $\T^3$:
$$\hat{\tau}(n_1,n_2,n_3) =\left(P_2(T^{n_1}f_1\otimes T^{n_2}f_2)\,,\, T^{n_3}f_3\right),$$
for fixed functions $f_1,f_2,f_3$ with zero mean, and $P_2$ is the Markov operator
(conditional mathematical expectation) corresponding to the self-joining $\nu$.

The support of such a measure $\tau$ lies in the subgroup $\{(t_1,t_2,t_3): t_1+t_2+t_3=0\}$,
moreover, $\pi_{i,j}\tau\ll \sigma\otimes\sigma$, where $\pi_{i,j}$ is the projection of the torus $\T_{1}\times\T_{2}\times \T_{3}$ onto the face $\T_{i}\times \T_{j}$, $1\leq i< j\leq 3$, and $\sigma$ is the measure of the maximal spectral type of the automorphism $T$.

Host established (the key point) that the projection $\pi\tau$ of such a measure $\tau$ onto $\T_1$ is the sum of a discrete and an absolute continuous measure, while $\pi\tau\ll \sigma$. But $\sigma$ is a continuous singular measure, so
$\pi\tau=0$ , $\tau=0$, and for  $f_1,f_2,f_3\perp Const$  we have
$$\int f_1\otimes f_2\otimes f_3\ d\nu = (P_2(T^{n_1}f_1\otimes T^{n_2}f_2)\, , \, T^{n_3}f_3)=0.$$
Thus, the self-joining $\nu$ is trivial.

\vspace{3mm}
\bf Theorem (B. Host \cite{Ho}). \it If spectral measure $\sigma$ of $T$ is singular,  then:

 \rm ($\ast$) \it any pairwise independent joining of automorphisms $R,S,T$ is trivial.\rm

\vspace{3mm} 
Exercise. Assume  additionally that the convolution power $\sigma^{\ast 2024}$ is absolutely continuous. Under this assumption find a short direct proof of  ($\ast$). 

The author proved the following assertion.

\vspace{3mm}
\bf Theorem (V. Ryzhikov \cite{Izv}). \it If  $T$ has $S_4$, then
 any pairwise independent joining of automorphisms $R,S,T$ is trivial.\rm

\vspace{3mm} 
So, there a  variation of del Junco - Rudolph
problem: \it  are there weakly mixing $R,S,T$ of zero entropy
possessing non-trivial  pairwise independent joining? \rm

\vspace{3mm}\it
Do $S_3$, $S_4$ imply weakly mixing and zero entropy for $\Z$-action?\rm

\vspace{3mm}
The one-dimensionality of the time gives this problem a peculiarity,
since for the $\Z^2$-actions counterexamples were discovered in advance by Ledrappier, see \cite{JR}.

\vspace{3mm}
\bf Bernoulli or Ledrappier. \rm  In \cite{Izv} there is the following observation. Let we have 
a non-trivial factor $F$  of ergodic $T\times T$  that is  independent with respect to the coordinate factors and generated by  a  set of the form 
$$(A\times A)\cup((X\setminus A)\times 
(X\setminus A).$$  Then $F$ is Bernoulli (Assertion 4.4. \cite{Izv}).   It was an exercise that  for $Z^2$-action the same factor  was Bernoulli or Ledrappier.

 \section{ Quantative deviation from mutiple mixing,  and local rank  of automorphisms}
We omit the stories related to Kalikow's pioneering result \cite{Ka} and the generalization of his theorem to finite-rank automorphisms \cite{93}. Let us  move on to the class of positive local rank transformations, where Rokhlin's problem is far from being solved.

An automorphism $T$ has\it local rank \rm
$\beta > 0$ if there is a sequence of towers 
 $$U_j=\sqcup_{k=1}^{h_j} T^kB_j, \ \ \mu(U_j)\to \beta, $$
  such that for    any $A\in\B$ the intersections $U_j\cap A$
 can be approximated by 
the union of some collections of floors in the tower  $U_j$. 
This means that
for some sequence  $S_j\subset \{0,1,\dots, h(j)-1\}$
(it depends on $A$) the  measures of 
$(U_j\cap A)\Delta\sqcup_{z\in S_j}T^zB$
tends to 0. 

Now we define \it  $(1+\varepsilon)$-mixing \rm -- a property that occupies 
an intermediate position between the property
of mixing  and the property of 2-fold mixing.
We set
$$
Dev(\varepsilon ,A,B,C) = 
 \lbrace (z,w)\in Q(\varepsilon ,h): | \mu (A\cap T^z B\cap
T^w C) - \mu (A)\mu (B)\mu (C)| > \varepsilon \rbrace,
$$
where
$$
Q(\varepsilon ,h) = \lbrace (z,w) \epsilon [0,h]\ : \
|z|,|w|,|z-w| > \varepsilon h \rbrace.
$$
Let
$$ 
d(h)=| Dev(\varepsilon ,A,B,C)|\, /\,  h.
$$

 In  \cite{92} it was  noted the finiteness of the deviations from multiple mixing:

\it for a mixing transformation $T$, sets $A,B,C$ a sequence $h_j\to\infty$,
$\eps>0$
$$ \limsup_j  {|\{i:\, \mu(A\cap T^{h_j}B\cap T^iC)\,>\, \eps\}|}\ <\infty.$$
  \rm

Thus, if  $T$ is mixing, then 
$d(h)$ is bounded.
If $T$ is 2-fold mixing  then $d(h)=0$
for large  $h$.  If for any  $A,B,C\in \B$ and $\varepsilon >0$
we have  $d(h) \rightarrow \, 0$, such  $T$  is called  \it
$(1+\varepsilon)$-mixing. \rm

From the results of \cite{97}, \cite{95}, \cite{97s} , and \cite{GHR}, \cite{KT},  
by hiding our  joinings, we obtain  the following assertion.  

\bf Theorem 5.1. \it  For a mixing automorphism $T$ with $\beta(T)>0$ the following 
properties {\rm (i),(ii),(iii)} are equivalent:

{\rm (i)}   $T$  is not 2024-fold mixing;

{\rm (ii)}    $T$  is   not  $(1+\eps)$-mixing, so  
 the number of deviations from 2-fold mixing is "maximal":  $\limsup_h d(h) > 0$
for some $A,B,C\in \B$ and $\varepsilon >0$;

{\rm (iii)}  there exist positive measure sets $A,B$  and  $m_i,n_i\to +\infty$ such that

\ \ \ \ \ \ $A\cap T^{m_i}A\cap T^{m_i+n_i}B =\varnothing.$

If {\rm (i)} holds, then    $T$  has Lebesgue spectrum of finite multiplicity.\rm

\large


\begin{thebibliography}{99}

\bibitem{Ro} V.A. Rokhlin, On endomorphisms of compact commutative groups, Izv. Akad. Nauk SSSR. Ser. Mat., 13:4 (1949), 329-340

\bibitem{L} F. Ledrappier, Un champ marcovien peut etre d'entropie null et melangeant, C.R. Acad. Sci. Paris Ser. A, 287 (1978), 561-563


\bibitem{An} H. Anzai, Ergodic skew product transformations on the torus. Osaka Math. J. 3 (1951), 83-89

\bibitem{R23} V.V. Ryzhikov, Generic extensions of ergodic systems, Sb. Math., 214:10 (2023), 1442-1457 


\bibitem{Leo} V.P. Leonov, Application of the characteristic functional and semi-invariants to the ergodic theory of stationary processes, Dokl. Akad. Nauk SSSR, 133:3 (1960), 523-526;
\bibitem{19}V.V. Ryzhikov, Weakly homoclinic groups of ergodic actions, Trans. Moscow Math. Soc., 80 (2019), 83-94


\bibitem{Ho}B. Host, Mixing of all orders and pairwise independent joinings of systems with singular spectrum, Israel J. Math., 76 (1991), 289-298


\bibitem{Ma}B. Marcus, The horocycle flow is mixing of all degrees. Inv. Math. 46 (1978), 201-209

\bibitem{91} V.V. Ryzhikov, Connection between the mixing properties of a flow and the isomorphicity of the transformations that compose it, Math. Notes, 49:6 (1991), 621-627

 \bibitem{Mo} S. Mozes, Mixing of all orders of Lie group actions. Invent. Math. 107 (1992), 235-241 

\bibitem{St}A.N. Starkov, On mixing of higher degrees of homogeneous flows. (Russian)Dokl. Ross. Akad. Nauk 333 (1993), 28-31 




\bibitem{Ra}M. Ratner, Horocycle flows, joinings and rigidity of products. Ann. Math. 118 (1983), 277-313 

\bibitem{Th}J.-P. Thouvenot, Some properties and applications of joinings in ergodic theory, Ergodic Theory and Its Connections with Harmonic Analysis, Proceedings of the Alexandria 1993 Conference, LMS Lecture Note Series, 205, eds. K. E. Petersen and I. A. Salama, 1995, 207-235

\bibitem{FK}B. Fayad, A. Kanigowski, Multiple mixing for a class of conservative surface flows, Invent. Math., 203:2 (2016), 555-614 

\bibitem{RT}
V.V. Ryzhikov, J.-P. Thouvenot, Disjointness, Divisibility, and Quasi-Simplicity of Measure-Preserving Actions, Funct. Anal. Appl., 40:3 (2006), 237-240 

\bibitem{00}  V.V. Ryzhikov, Rokhlin's multiple mixing problem in the class of positive local rank actions, Funct. Anal. Appl., 34:1 (2000), 73-75

\bibitem{Ka}S.A. Kalikow, Twofold mixing implies threefold mixing for rank one transformations. Ergod. Theory Dynam. Syst. 4 (1984), 237-259

\bibitem{93} V.V. Ryzhikov, Joinings and multiple mixing of finite rank actions. Funct. Anal. Appl. 27 (1993), No. 2, 128-140

\bibitem{KG} J.L. King, On M. Gordin's homoclinic question, Internat. Math. Res. Notices, 1997, no. 5, 203-212


\bibitem{JR}A. del Junco, D. Rudolph, On ergodic action whose self-joinings are graphs, Ergodic Theory Dynam. Systems, 7 (1987), 531-557

\bibitem{Ki}J. King, Ergodic properties where order 4 implies infinite order. Israel J. Math. 80, 1-2 (1992), 65-86



\bibitem{97} V.V. Ryzhikov, Intertwinings of tensor products, and the stochastic centralizer of dynamical systems. Sb. Math. 188, No. 2 (1997), 237-263


\bibitem{92}V.V. Ryzhikov, Mixing, rank, and minimal self-joining of actions with an invariant measure, Russian Acad. Sci. Sb. Math., 75:2 (1993), 405-427  

\bibitem{GHR}E. Glasner, B. Host, and D. Rudolph, Simple systems and their higher order self-joinings. Israel J. Math. 78 (1992), 131-142

\bibitem{95} V.V. Ryzhikov, Multiple Mixing and Local Rank of Dynamical Systems. Funct. Anal. Appl., 29:2 (1995), 143-145


\bibitem{97s}V.V. Ryzhikov, Around simple dynamical systems. Induced joinings and multiple mixing. J. Dyn. Control Syst. 3 (1997), No. 1, 111-127

\bibitem{KT}J.L. King and J.-P. Thouvenot, A canonical structure theorem for finite joiningrank maps, J. Analyse Math. 56 (1991), 211-230

\bibitem{Ti}S. V. Tikhonov, On the absence of multiple mixing and on the centralizer of measure-preserving actions, Mat. Notes, 97:4 (2015), 636-640

\bibitem{Izv} V.V. Ryzhikov, Joinings, intertwining operators, factors, and mixing properties of dynamical systems, Russian Acad. Sci. Izv. Math., 42:1 (1994), 91-114  
\end{thebibliography}
\end{document}